\documentclass{article}
\usepackage{amsmath, amssymb, amsfonts}
\newtheorem{thm}{Theorem}

\newtheorem{cor}{Corollary}
\newtheorem{lemma}{Lemma}

\newenvironment{proof}{\medskip 
\noindent {\bf Proof.}}{\hfill \rule{.5em}{1em}\mbox{}\bigskip}
\def\Bbb{\mathbb} 
\def\bcp{{\mathbb C \mathbb P}}
\newtheorem{main}{Theorem}

\begin{document}

\title{On Einstein Manifolds of\\
Positive Sectional Curvature}

\author{
Matthew J. Gursky\thanks{Supported 
in part by  NSF grant DMS-9623048.}
\\ Indiana University
\\ and \\
Claude LeBrun\thanks{Supported 
in part by  NSF grant DMS-9802722.}\\
SUNY Stony Brook
}
\date{July 8, 1998}

\maketitle

 \begin{abstract}
Let $(M,g)$ be a compact oriented $4$-dimensional 
Einstein manifold. If $M$ has positive intersection form
and $g$ has non-negative sectional curvature, we show that,
up to rescaling and isometry, $(M,g)$ is $\bcp_{2}$, 
with its standard  Fubini-Study metric. 
 \end{abstract}

\section{Introduction}

A  Riemannian manifold $(M,g)$ is said to be {\em Einstein}
if it has constant Ricci curvature --- i.e. if
its  Ricci tensor $r$ is a constant multiple of 
the metric: 
\begin{equation}
	r=\lambda g .
	\label{eins}
\end{equation}
If $g$ is complete and $\lambda > 0$, Myers' theorem 
\cite{myers} then 
tells us that $M$ is compact, and has finite fundamental group.

The simplest examples of compact Einstein manifolds with 
positive Ricci curvature 
($\lambda > 0$) are provided by the irreducible symmetric
spaces of compact type. In dimension $4$, this observation 
yields  exactly
two orientable examples: $S^{4}=SO(5)/SO(4)$ and 
$\bcp_{2}=SU(3)/U(2)$, both of which actually
have positive {\em sectional} curvature. A slight generalization 
would be  to allow for {\em reducible} symmetric spaces; in dimension $4$,
this gives us the additional oriented examples of $S^{2}\times S^{2}=
SO(4)/[SO(2)\times SO(2)]$ and its quotient by the simultaneous
antipodal map on both factors. The latter examples have non-negative
sectional curvature, although some of their sectional curvatures
are actually zero. 

While \cite{tian}  there certainly are other compact $4$-dimensional  
Einstein manifolds with $\lambda > 0$, none are known which have 
 non-negative sectional
curvature. One might hope that  this is not merely
accidental. In this direction, we are able to offer the following partial 
result:

\begin{main}\label{A}
Let $M$ be a smooth compact oriented 4-manifold with
(strictly) positive intersection form, and 
suppose that $g$ is an Einstein metric on $M$
which has non-negative sectional curvature. 
Then  $(M,g)$ is homothetically isometric
 to  $\bcp_{2}$, equipped with its standard Fubini-Study metric. 
\end{main}

To clarify the statement, let us  
recall that one can always find  bases for the de Rham 
cohomology $H^{2}(M, {\Bbb R})$ in  which the intersection 
pairing 
\begin{eqnarray*}
H^{2}(M, {\Bbb R})\times H^{2}(M, {\Bbb R})	
 & \longrightarrow & ~~~~ {\Bbb R}  \\
	( ~ [\varphi ] ~ , ~ [\psi ] ~) ~~~~~ & 
	\mapsto  & \int_{M}\varphi \wedge \psi 
\end{eqnarray*}
is represented by a diagonal matrix  
	   
	    $$\left[ 
	  \begin{array}{rl}  \underbrace{
	 \begin{array}{ccc}
	   		 1 &  &	  \\
	   		  &	\ddots &   \\
	   		  &	 & 1
	   	  \end{array}}_{b_{+}(M)}
	   	   & 
	   		   \\ 
	  
	 {\scriptstyle b_{-}(M)}  \!
	    \left\{\begin{array}{r}
	    \\
	    \\
	    \\
	    \\
	    \end{array}
	    \right. \! \! \! \! \! \! \! \! \! \! \! \! \! \! \! 
	    &\begin{array}{ccc}
	   		 -1	&  &   \\
	   		  &	\ddots &   \\
	   		  &	 & -1
	   	  \end{array}
	    \end{array} 
	    \right] 
	    $$
and, in these terms, our  first hypothesis  
stipulates that $b_{-}(M)=0$ and that $b_{+}(M) \neq 0$.  
The {\em Fubini-Study metric} is the unique $U(2)$-invariant metric on
$\bcp_{2}= SU(3)/U(2)$ with total volume $\pi^{2}/2$; it is 
Einstein, and  has
 sectional curvatures $K(P)\in [1,4]$. 
By {\em homothetically isometric}, we mean isometric
after rescaling; in other words, the  
theorem concludes by  asserting the existence of 
 a diffeomorphism $\Phi : M \to \bcp_{2}$ such that 
$g=\Phi^{*}cg_{0}$ for some  some positive constant $c$.

Theorem \ref{A} is actually a consequence of the following,
more general result:

\begin{main}\label{B}
Let $(M,g)$ be a smooth compact oriented Einstein 4-manifold
with  non-negative sectional curvature. Assume, moreover, that
$g$ is neither self-dual nor anti-self-dual.  Then  
the Euler characteristic $\chi$  and the signature $\tau$ of $M$ 
satisfy 
$$9\geq \chi  > \frac{15}{4} |\tau |.$$
\end{main}
Here $\tau (M) := b_{+}(M)-b_{-}(M)$, whereas 
$\chi (M)= 2 +  b_{+}(M)+b_{-}(M)$ if $M$ has finite fundamental group. 
 Thus, for example,  while Tian \cite{tian} has shown that the manifolds 
 $\bcp_{2}\#k\overline{\bcp}_{2}$, $3\leq k  \leq  8$, 
 admit Einstein 
 metrics with positive Ricci curvature, 
these spaces definitely {\em do not} admit Einstein metrics of positive
{\em sectional} curvature. 
 
 Theorem \ref{B} was directly inspired by the following 
 result of Hitchin \cite{hit}: if a compact oriented 4-manifold
 $M$ admits an Einstein metric $g$ of positive sectional curvature,
 then the Euler characteristic $\chi$ and signature $\tau$ of $M$
 must satisfy
 $$\chi \geq \left(\frac{3}{2}\right)^{3/2}|\tau |.$$
 By freely quoting the subsequent topological results of Freedman 
 \cite{f}, 
 this actually contains enough information to conclude that an Einstein 
 4-manifold of
 positive
 intersection form and  positive sectional curvature
  must be {\em homeomorphic} to either
 $\bcp_{2}$ or $\bcp_{2}\# \bcp_{2}$. Even on the crude level of 
 homeomorphism type, however,  Theorem \ref{A}  represents 
 a 2-fold improvement over Hitchin's result. This
  reflects the  fact
 that Hitchin's coefficient of 
  $\left(\frac{3}{2}\right)^{3/2}= 1.837\ldots $ is less than half 
 the coefficient of $\frac{15}{4}= 3.75$ appearing in Theorem \ref{B}.

 Let us conclude this introduction with some
 general remarks regarding the class of Riemannian manifolds
 under consideration. First of all, a celebrated result of 
 Synge \cite{synge} asserts that any compact, orientable, 
  even-dimensional Riemannian 
 manifold of positive sectional curvature is necessarily
 simply connected. On the other hand, we have already seen,
 by the example of $(S^{2}\times S^{2})/{\Bbb Z}_{2}$, that
 no such result holds when
  the sectional curvature is 
 merely assumed to be non-negative, even if the Ricci curvature is
 positive. Nonetheless, Theorem \ref{B}
 {\em does} tell us that a compact orientable Einstein 4-manifold
 of non-negative sectional curvature and {\em non-zero signature}
 must be simply connected. Indeed, if $|\tau | \geq 1$, the 
inequality $\chi > \frac{15}{4}|\tau |$ and the observation that 
$\chi \equiv \tau \bmod 2$ together then  guarantee that 
$\chi \geq 5$. But if such a manifold were not simply connected,
its universal cover would then violate the inequality $9 \geq \chi$. 

Finally, let us observe that there are, up to diffeomorphism,
only finitely many compact 4-manifolds
with Einstein metrics of non-negative Ricci curvature. 
The flat 4-manifolds, of course, nominally 
 form a subclass of the the manifolds under
discussion, 
but Bieberbach's theorem \cite{bie} in any case tells us that
there are finitely many diffeomorphism types of these. 
For the others, which are our real  concern here, 
the Ricci curvature must be
positive, and we may thus rescale  the metric so that, for example,
 $r= 3g$. The definition of the Ricci curvature then tells us that  
 the sectional curvatures  all satisfy $0\leq K(P)\leq 3$.
 On the other hand, Myers' theorem \cite{myers} 
 predicts that the diameter is 
  $\leq \pi$. Moreover,  the 4-dimensional Gauss-Bonnet theorem
  \cite{bes,hit} 
and   our upper bound on curvature
imply the volume is 
$\geq 8\pi^{2}/15$. With such bounds, 
 Cheeger's finiteness theorem \cite{cheeger} then predicts that 
the given class of manifolds is precompact in the 
$C^{\alpha}$ topology\footnote{Indeed, the  ellipticity of 
 the  
 Einstein equations guarantees \cite{anderson}  that
  the class of manifolds in question is 
actually compact in the $C^{\infty}$ topology.},
 and  therefore consists of  finitely many 
 diffeomorphism classes.  
 Unfortunately, however, such arguments by no means 
  predict the actual number of diffeotypes.
By contrast, Theorem \ref{B} and Freedman's classification 
\cite{f} 
tell us that there are at most twelve {\em homeotypes}
of  simply connected compact  Einstein 4-manifolds with non-negative
sectional curvature.

 \section{The Curvature of 4-Manifolds}
 \label{local}
 
 We begin by
 recalling  that the 
rank-$6$ bundle of $2$-forms $\Lambda^{2}$ on 
 an oriented  Riemannian 4-manifold $(M^{4},g)$
has an invariant decomposition 
\begin{equation} 
\Lambda^2=\Lambda^+\oplus \Lambda^-
\label{deco} 
\end{equation}
as the sum of two rank-$3$ vector  bundles.
Here $\Lambda^{\pm}$ are by definition the 
eigenspaces of the  Hodge duality
operator
$$\star: \Lambda^2 \to \Lambda^2,$$
corresponding respectively  to the  eigenvalue $\pm 1$. 
Sections of $\Lambda^{+}$ are called {\em self-dual}
2-forms, whereas sections of $\Lambda^{-}$ are called
 {\em anti-self-dual} 2-forms. 
But since the curvature tensor of $g$ may be 
thought of as a map
${\cal R} : \Lambda^2 \to \Lambda^2$,  
 (\ref{deco}) gives us a decomposition \cite{st}  of the curvature into 
  primitive pieces 
\begin{equation}
\label{curv}
{\cal R}=
\left(
\mbox{
\begin{tabular}{c|c}
&\\
$W^++\frac{s}{12}$&$\stackrel{\circ}{r}$\\ &\\
\cline{1-2}&\\
$\stackrel{\circ}{r}$ & $W^-+\frac{s}{12}$\\&\\
\end{tabular}
} \right) , 
\end{equation}
where the {\em  self-dual} and 
{\em anti-self-dual Weyl curvatures}  
$W^\pm$  
are trace-free   as endomorphisms of $\Lambda^{\pm}$.
The scalar curvature  $s$ is understood here
  to act by scalar multiplication. On the other hand, 
  $\stackrel{\circ}{r}$ represents the 
  trace-free Ricci curvature
  $r-\frac{s}{4}g$, and so vanishes iff $g$ is Einstein.
  
This last fact has a simple but crucial consequence.

\begin{lemma} Let $(M,g)$ be an oriented Einstein $4$-manifold. 
If the sectional curvature of $g$ is non-negative, then 
\begin{equation}
\frac{s}{\sqrt{6}} \ge |W^+| + |W^-|
\label{2.3}
\end{equation}
at each point of $M$. \label{K} 
\end{lemma}

\begin{proof}
Every 2-form $\varphi$ on $M$ can be uniquely written as
$\varphi = \varphi^{+} + \varphi^{-}$,
where $\varphi^{\pm}\in \Lambda^{\pm}$. Now a 
2-form is expressible as a simple wedge product of 
1-forms iff $\varphi \wedge \varphi =0$.
But this condition can be rewritten as 
$|\varphi^{+}|^{2}-|\varphi^{-}|^{2}=0$. 
Thus the sectional curvature of $g$ is non-negative 
iff the curvature operator ${\cal R} : \Lambda^{2}\to \Lambda^{2}$
satisfies
$$\langle \varphi^{+} + \varphi^{-},
 {\cal R}(\varphi^{+} + \varphi^{-})\rangle \geq 0$$
for all {\em unit-length} self-dual 2-forms  $\varphi^{+}$
and all {\em unit-length} anti-self-dual 2-forms  $\varphi^{-}$. 
But for an Einstein manifold, (\ref{curv}) tells us that this can be 
rewritten as 
\begin{equation}
\frac{s}{6} + \lambda_{+} + \lambda _- \ge 0
\label{2.2}
\end{equation}
where, for each $x\in M$, 
 $\lambda _\pm(x)\leq 0$ is by definition the smallest eigenvalue of 
the trace-free endomorphism 
$W^{\pm}_{x}: \Lambda^{\pm}_{x}\to \Lambda^{\pm}_{x}$.

The claim will thus follow  immediately from (\ref{2.2}) 
if we can show that
\[
|\lambda _\pm| \geq  \frac{1}{\sqrt{6}} |W^\pm |.
\]
To see, this, let 
 $\lambda_{+}\leq \mu_+ \leq \nu_+$ be the  eigenvalues of $W^+$.
Thus 
\[
|W^+|^2 = \lambda ^2_+ + \mu ^2_+ + \nu ^2_+.
\]
But since $W^+$ is
trace-free, $\lambda _+ + \mu _+ + \nu _+ =0$,
and hence 
\begin{eqnarray*}
|W^+|^2	 & = &  \lambda ^2_+ + \mu ^2_+ + \nu ^2_+
+ (\lambda _+ - \mu _+ - \nu _+)(\lambda _+ + \mu _+ + \nu _+)  \\
	 &  = &  
2\left[ \lambda_{+}^{2}-\mu_{+}\nu_{+}
\right].
\end{eqnarray*}
If $\mu_{+}\geq 0$, this last expression is less than $2|\lambda_{+}|^{2}$.
 Otherwise, $\lambda_{+}\leq \mu_{+} < 0$, 
$0 < \nu_{+}  \leq 2 |\lambda_{+}|$, and 
hence 
$
|W^+|^{2}
\leq  6|\lambda _+|^{2}
$.
Thus 
$|\lambda_{+}|\geq \frac{1}{\sqrt{6}}|W^{+}|$.
Since $|\lambda_{-}|\geq \frac{1}{\sqrt{6}}|W^{-}|$
by the same argument, we are done. 
\end{proof}

The curvatures $W^{\pm}$, $\stackrel{\circ}{r}$, and $s$ correspond
to 
 different irreducible representation of 
$SO(4)$, so the only
 invariant quadratic polynomials in the curvature 
 of an oriented $4$-manifold are 
 linear combinations of $s^{2}$, $|\stackrel{\circ}{r}|^{2}$,
$|W_{+}|^{2}$ and $|W_{-}|^{2}$. This observation can be 
applied, in particular, to simplify the integrands
\cite{bes,hit} of the
4-dimensional  Chern-Gauss-Bonnet 
\begin{equation}
\chi (M)= \frac{1}{8\pi^2}\int_M \left[|W_+|^2
+ |W_-|^2 + \frac{s^2}{24} -\frac{|\stackrel{\circ}{r}|^2}{2}
\right]d\mu 
		\label{gb}
\end{equation}
and Hirzebruch signature   
 \begin{equation}
 	\tau (M)= \frac{1}{12\pi^2}\int_M \left[|W_+|^2
- |W_-|^2  \right]d\mu 
 	\label{sig}
 \end{equation}
  formul{\ae}.
Here the curvatures, norms $|\cdot |$, and
volume form $d\mu$ are,  of course, those of any given  
Riemannian metric $g$ on $M$. 

Applying  Lemma \ref{K} now gives us 
an  elementary but useful result:

\begin{lemma}\label{chi}
Let $(M,g)$ be a compact  4-dimensional Einstein manifold
of non-negative sectional curvature. If $g$  is not flat, then
\begin{equation}
\chi (M) < \frac{5}{8\pi^{2}} \int_{M}\frac{s^{2}_{g}}{24}	d\mu_{g}.
	\label{euler}
\end{equation}
\end{lemma}
\begin{proof}
 Because $g$ is not flat, and  the sectional curvature is
non-negative, our Einstein metric $g$ must have positive 
scalar curvature, and hence positive Ricci curvature.  Myers' 
Theorem \cite{myers} thus forces  
$M$ to have finite fundamental group, so that,  in particular, 
$b_{1}(M)=0$. 
 
By passing to a double cover if necessary, we may assume 
that $M$ is orientable. Let us choose to  orient $M$ so that 
$\tau (M) \geq 0$.

 Lemma \ref{K} now  tells us that 
$$|W^{+}_{g}|^{2}+ |W^{-}_{g}|^{2} \leq \left(
|W^{+}_{g}| + |W^{-}_{g}|
\right)^{2} \leq \frac{s_{g}^{2}}{6},$$
so that 
$$
\chi (M) = \frac{1}{8\pi^{2}}\int_{M}\left[ 
|W^{+}_{g}|^{2}+|W^{-}_{g}|^{2}+\frac{s^{2}_{g}}{24}
\right] d\mu_{g} \leq 5\cdot \frac{1}{8\pi^{2}}\int_{M}\frac{s_{g}^{2}}{24}
 .$$
 
 If equality were to  hold,  we would have 
 $|W^{+}_{g}|~|W^{-}_{g}|\equiv 0$ and $|W^{+}_{g}| + 
 |W^{-}_{g}| \equiv \frac{s}{\sqrt{6}}$. But any Einstein metric $g$ is	
 real-analytic \cite{deka}
  in harmonic coordinates. With our orientation conventions, 
  we must therefore have
   $|W^{-}_{g}|\equiv 0$. Our Gauss-Bonnet 
    formul{\ae} (\ref{gb}--\ref{sig}) then tell 
   us that $\chi (M)= \frac{15}{8}\tau (M).$ However, 
   the Weitzenb\"ock formula for harmonic 2-forms implies \cite{bourg}
   that a self-dual 4-manifold with $s > 0$ has $b_{-}=0$.
   Since we also have $b_{1}(M)=0$, it follows that 
   $\chi (M) = 2 + \tau (M)$. Solving for the signature,
   we find that  $\tau (M) = \frac{16}{7}$.
   As  
   this is  of course
   a contradiction, it follows that  the inequality is always strict. 
\end{proof}

This has an important consequence: 

\begin{lemma}
\label{bshp}
Let $(M,g)$ be a compact  4-dimensional Einstein manifold
of non-negative sectional curvature. Then $\chi (M) \leq 9$. 
\end{lemma}
\begin{proof}
We may  assume that $g$ has positive Ricci 
curvature,  since otherwise the Euler characteristic 
would vanish.  
By rescaling, we can thus 
  arrange for our Einstein metric to have Ricci 
tensor $r=3g$.   
Bishop's inequality 
\cite{bishop} then asserts that the total volume of $(M,g)$
is less than or equal to that of  the 4-sphere with its 
standard metric $g_1$. 
Since both $g$ and $g_1$ have $s=12$,  Lemma \ref{chi} now asserts that  
$$\chi (M) < \frac{5}{8\pi^{2}} \int_{M}\frac{s^{2}_{g}}{24}	d\mu_{g}
\leq  \frac{5}{8\pi^{2}}\int_{S^{4}}\frac{s_{g_1}^{2}}{24}d\mu_{g_1} 
= 5\chi (S^{4}) = 10. 
$$
Since the Euler characteristic is an integer, it follows  that 
$\chi (M) \leq 9$.
\end{proof}

\section{$L^{2}$ Curvature  Estimates}
\label{global}

The key observations of \S \ref{local} were  basically 
 point-wise in  character. We now turn to some results 
of a fundamentally global nature, 
 beginning with a simplified 
  proof of a surprising fact   discovered in 
 \cite{G1}. 

\begin{lemma}  
  Suppose $(M,g)$ is a compact oriented 
 Einstein $4$-manifold of positive scalar
curvature.  Then either $W^+\equiv 0$, or else 
there is a smooth, conformally related  metric $\hat{g}=u^{2}g$
such that 
$$\int_{M}\left[ s_{\hat{g}}
- 2\sqrt{6} | W^{+}_{\hat{g}}|_{\hat{g}}\right] d\mu_{\hat{g}} \leq 0
 .$$
Moreover,  one can either  arrange for the inequality to be 
strict,  or for the metric $\hat{g}$ to be locally 
K\"ahler. 
\label{weyl} 
\end{lemma}

\begin{proof}
For each metric $\hat{g}$ on our oriented $4$-manifold $M$,
let us consider the quantity ${\frak S}_{\hat{g}}$ 
defined by 
$$
{\frak S}_{\hat{g}} = s_{\hat{g}} - 2\sqrt{6} |W^+_{\hat{g}}|_{\hat{g}} .
$$
Under conformal rescaling,  this curvature function behaves very much
 like the usual 
 scalar curvature $s$. 
 Indeed, 
 if $\hat{g}=u^2g$, where 
  $u$ is a smooth positive function,  
 then 
\begin{equation}
{\frak S}_{\hat{g}}= 
u^{-3}\diamondsuit u ,
\label{transf} 
\end{equation}
where,
in terms of  
the   (positive) Laplace-Beltrami operator
$\Delta =d^*d=-\mbox{div } \mbox{grad}$, 
the linear  elliptic operator $\diamondsuit = \diamondsuit_{g}$ 
is defined by 
 $$
\diamondsuit 
=
6\Delta_g+  {\frak S}_g  .
$$
Since $d\mu_{\hat{g}}= u^{4}d\mu_{g}$, we thus have 
$$\int_{M}\left[ s_{\hat{g}}
- 2\sqrt{6} | W^{+}_{\hat{g}}|_{\hat{g}}\right] d\mu_{\hat{g}}=
\int_{M}	{\frak S}_{\hat{g}}	d \mu_{\hat{g}}
= \int_{M} (u\diamondsuit u)d\mu_{g}.
$$

The above generalities apply to any conformally related pair of 
metrics. But 
in the present case, the given metric $g$ is assumed to be Einstein.
The second Bianchi identity therefore tells us that 
 its self-dual Weyl curvature 
 is harmonic,  in the sense that 
\begin{equation}
\nabla^{a}W^{+}_{abcd}=0. 
\label{div}
\end{equation}
In spinor terms, this says that  $\nabla W^{+}\in {\Bbb S}_{-}\otimes
\bigodot^{5}{\Bbb S}_{+}$. Now suppose   
 $U\in \bigodot^{4} {\Bbb S}_{+}$ and 
 $v\in {\Bbb S}_{-}\otimes {\Bbb S}_{+} ={\Bbb C}\otimes TM$
 are real elements,
and let $(v\otimes  U)^{\|}$ denote the orthogonal projection 
of $v\otimes  U$ to ${\Bbb S}_{-}\otimes
\bigodot^{5}{\Bbb S}_{+}$. Using the notational conventions of 
\cite{pr}, we then have  
\begin{eqnarray*}
	\left[ (v\otimes  U)^{\|}\right]_{A'ABCDE} & = &  v_{A'(A}U_{BCDE)}  \\
	 & = & \frac{1}{5}\left[ v_{A'A}U_{BCDE}+
	 v_{A'B}U_{ACDE}\right. \\& & 
	\left. \hphantom{{\varepsilon_{C}}^{A}}+
v_{A'C}U_{ABDE}+v_{A'D}U_{ABCE}+v_{A'E}U_{ABCD} \right] ,
\end{eqnarray*}
so that 
\begin{eqnarray*}
|(v\otimes  U)^{\|}|^{2}	 & = & 
v^{A'A}U^{BCDE}v_{A'(A}U_{BCDE)}
  \\
	 & = & \frac{1}{5}v^{A'A}U^{BCDE}\left[ v_{A'A}U_{BCDE}+
	 v_{A'B}U_{ACDE}+v_{A'C}U_{ABDE}\right.
	 \\& &  \left. \hphantom{ABCDEFGHIJK}+v_{A'D}U_{ABCE}+v_{A'E}U_{ABCD}
	   \right] 
  \\ &=&
   \frac{|v|^{2}}{5}U^{BCDE}\left[ U_{BCDE}+\frac{1}{2}{\varepsilon_{B}}^{A}
  U_{ACDE}+
\frac{1}{2}{\varepsilon_{C}}^{A}U_{ABDE}+
\right. \\& & \left. \hphantom{ABCDEFGHI}+
  \frac{1}{2}{\varepsilon_{D}}^{A}U_{ABCE}+
  \frac{1}{2}{\varepsilon_{E}}^{A}U_{ABCD} \right] \\
	 & = &  \frac{3}{5}|v|^{2}|U|^{2}. 
\end{eqnarray*}
The Cauchy-Schwarz inequality therefore predicts that 
$$\langle v\otimes  U , \nabla W^{+}\rangle \leq \sqrt{\frac{3}{5}}
|v|~|U|~|\nabla W^{+}|. $$
Away from the zeroes of $W^{+}$,
setting $U=W^{+}$ thus yields  
$$|W^{+}| ~\nabla_{v} |W^{+}| = 
 \langle v\otimes  W^{+} , \nabla W^{+}\rangle
 \leq  \sqrt{\frac{3}{5}} ~ |v|~|W^{+}|~|\nabla W^{+}|,$$
giving us the Kato inequality
\begin{equation}
|\nabla W^{+}| \geq \sqrt{\frac{5}{3}}~ \left| \nabla |W^{+}| \right|. 
\label{kato}
\end{equation}

On the other hand, 
Derdzi\'nski \cite{derdz} observed that equation
(\ref{div}) also implies  the
Weitzenb\"ock formula
\begin{equation}
0=\frac{1}{2} \Delta |W^+ |^2  + |\nabla W^+|^2 + \frac{s}{2}
|W^+|^2 - 18 \det W^+
\label{1.3}
\end{equation}
where $\Delta$ is again the  positive Laplacian and 
$\det W^+$ is the determinant of the bundle endomorphism
$W^+:\Lambda^+ \to \Lambda^+$; cf. \cite[equation (6.8.40)]{pr}.   
In conjunction with 
 the (sharp) algebraic inequality $3 \sqrt{6}
\det W^+ \leq  |W^+|^3$, equations  (\ref{kato}--\ref{1.3}) 
imply that the non-negative function $u_{0}= |W^+|^{1/3}$ satisfies 
\begin{equation}
0\geq \diamondsuit u_{0} 
\label{1.4a}
\end{equation} 
in the classical sense, 
except at the locus $u_{0}=0$, where it presumably fails to be smooth. 
Now, for each $\epsilon > 0$,
let $f_{\epsilon }: [0,\infty ) \to (0,\infty )$ be a smooth
positive  
function which is constant on $[0,\epsilon /2 ]$, satisfies
$f_{\epsilon }(x)=x$ for $x >  \epsilon$, and has non-negative
 second derivative everywhere. We may then consider the 
 smooth positive function  $u_{\epsilon} = f_{\epsilon 
}\circ u_{0}$, and the metric $g_{\epsilon}= u_{\epsilon}^{2}g$. 
 Let $M_{\epsilon}$ be the 
set where $u_{0} < \epsilon$. Then 
\begin{eqnarray*} \int_{M}{\frak S}_{g_{\epsilon}}
d\mu_{g_{\epsilon}}
	 & = &  \int_{M} (u_{\epsilon} \diamondsuit
	  u_{\epsilon} ) d\mu_{g} \\
	 &  \leq & C \epsilon^{2}
\mbox{Vol} (M_{\epsilon})
+ \int_{M-M_{\epsilon}} 
(u_{0}  \diamondsuit u_{0}  ) d\mu_{g} ,
\end{eqnarray*}
where $C$ is any positive  upper bound for ${\frak S}_{g}$. 

Now assume that $W^{+}\not\equiv 0$. 
Since $g$ is real-analytic \cite{deka} in harmonic coordinates, 
so is $u_{0}^{6}=|W^{+}|^{2}$, and hence 
  $\mbox{Vol} (M_{\epsilon}) \longrightarrow 0$
  as $\epsilon \to 0$. Thus 
  $\int_{M}	{\frak S}_{g_{\epsilon}}	d \mu_{g_{\epsilon}}$
 is  negative for small $\epsilon$ unless equality holds in (\ref{1.4a}). 
 
 In the latter case, however, we have
  $|W^{+}|^{3}\equiv 3\sqrt{6}\det W^{+}$, 
  and hence  $W^{+}$  has at most
 2 distinct eigenvalues at each point. Derdzi\'nski's
theorem \cite{derdz} thus asserts that $W^{+}\neq 0$, and 
that $\hat{g}= u_{0}^{2}g$ is locally K\"ahler --- i.e. becomes
 K\"ahler after possibly
pulling back to a double cover of $M$. 
\end{proof}

This implies  \cite{G1}  a remarkable ``gap theorem'' for $W^{+}$:  

\begin{thm}  
 Let $(M,g)$ be a compact oriented 
 Einstein $4$-manifold with $s >0$ and  $W^+\not\equiv 0$.
 Then  
$$ 
\int_{M} |W^{+}_{g}|^{2}_{g}d\mu_{g} \geq 
\int_{M}\frac{s^{2}_{g}}{24}d\mu_{g}
 ,$$
 with equality iff $\nabla W^{+}\equiv 0$. \label{gap} 
\end{thm}
  \begin{proof}
  A fundamental result of Obata \cite{obata}
 implies that  any   Einstein metric 
  is   a Yamabe minimizer;  moreover, such a  metric is always the 
  {\em unique}  Yamabe 
  minimizer, modulo homotheties and --- on the round 
   sphere --- global conformal transformations.
  Thus, if $\hat{g}=u^{2}g$ is any conformal rescaling of our
  Einstein metric $g$, we have  
  $$
  \frac{
  \int_{M} s_{{g}}  d\mu_{{g}} 
  }{
  \sqrt{\int_{M}
 d\mu_{{g}}}
 }
\leq 
  \frac{
  \int_{M} s_{\hat{g}}  d\mu_{\hat{g}} 
  }{
  \sqrt{\int_{M}
 d\mu_{\hat{g}}}
 } .
  $$ 
  However, assuming that $W^{+}\not\equiv 0$,
 Lemma \ref{weyl} tells us that $u$ can be chosen so that 
  \begin{eqnarray*}
  	\int s_{\hat{g}} d\mu _{\hat{g}}  & \leq &
  	 2 \sqrt{6} \int |W^+_{\hat{g}}|d\mu _{\hat{g}}  \\
  	 & \leq & \left( 24\int|W^+_{\hat{g}}|^2 d
\mu_{\hat{g}}\right)^{\frac{1}{2}} \left(\int d\mu
_{\hat{g}}\right)^{\frac{1}{2}}.
  \end{eqnarray*}
  Since $s_{g}$ is constant, and because 
  the $L^{2}$ norm of $W^{+}$ is conformally invariant, it therefore
  follows that 
  \begin{eqnarray*}
  	\left(\int_{M}s^{2}_{g}d\mu_{g}
  	\right)^{1/2} & = &  
  	\frac{
  \int_{M} s_{{g}}  d\mu_{{g}} 
  }{
  \sqrt{
  \int_{M}
 d\mu_{{g}}
 }
 } \\ 
 & \leq  &   \frac{
  \int_{M} s_{\hat{g}}  d\mu_{\hat{g}} 
  }{
  \sqrt{
  \int_{M}
 d\mu_{\hat{g}}
 }
 }  \\ & \leq  &  \left( 24\int|W^+_{\hat{g}}|^2 d
\mu_{\hat{g}}\right)^{\frac{1}{2}}
  	 \\  & = &  \left(24 \int |W^+_g|^2 d\mu _g
  	 \right)^{\frac{1}{2}} .
  \end{eqnarray*} 
 Moreover, equality can occur only if 
 $\hat{g}$ is both locally K\"ahler and isometric to a constant times  
  $g$. The latter, of course, happen iff
 $g$ is itself locally K\"ahler. But since $s\neq 0$ is constant and
 $W^{+}\not\equiv 0$,  the latter is equivalent to requiring that 
 $\nabla W^{+}\equiv 0$. \end{proof}

Reading this  in the mirror, we have: 

\begin{cor}\label{gapkid}
Let $(M,g)$ be an   oriented  compact Einstein $4$-manifold
with $s > 0$ and  $W^{-}\not\equiv 0$. Then 
 \begin{eqnarray}
 \int_{M} |W^-_g|^2 d \mu_g 	 &  \geq  & 
 \int \frac{s^2_g}{24} d\mu_g  ~ \mbox{\it  , and} 
 	\label{i}  \\
 	\frac{2\chi - 3\tau }{3} (M) & \geq  & \frac{1}{4\pi^{2}}
 	\int_{M} \frac{s^2_g}{24} d\mu_{g} ~ . 
 	\label{ii}
 \end{eqnarray}
Moreover, both these inequalities are strict unless $\nabla W^{-}\equiv 0$.
\end{cor}
\begin{proof}
  Reversing the orientation of $M$ interchanges $W^{+}$ and $W^{-}$.
  Applying this observation to  Theorem   \ref{gap} immediately yields 
  (\ref{i}).  
  But this and the 
   Gauss-Bonnet-type 
   formul{\ae} (\ref{gb}--\ref{sig}) then 
   tell us that 
  \begin{eqnarray*}
   	(2 \chi -3\tau)(M) & = & \frac{1}{4\pi^{2}}\int_{M}
\left[ 2 |W^-_g|^2 +\frac{s^{2}}{24} \right] d \mu_g  \\
   	 & \geq  & \frac{3}{4\pi^{2}}
 	\int_{M} \frac{s^2_g}{24} d\mu_{g} ~,
   \end{eqnarray*}
  thus proving (\ref{ii}). 
 \end{proof}

\section{The Main Theorems}

Combining the estimates of \S\S \ref{local}--\ref{global}
now  allows us to prove our main inequality:

  \setcounter{main}{1}
\begin{main}
Let $(M,g)$ be a smooth compact oriented Einstein 4-manifold
with  non-negative sectional curvature. Assume, moreover, that
$g$ is neither self-dual nor anti-self-dual.  Then  
the Euler characteristic $\chi$  and the signature $\tau$ of $M$ 
satisfy 
$$9\geq \chi  > \frac{15}{4} |\tau |.$$
\end{main}
\begin{proof}
 Combining (\ref{euler}) and (\ref{ii}), we have 
 $$\frac{2}{3}\chi - \tau \geq  \frac{1}{4\pi^{2}}
 \int \frac{s^{2}}{24}d\mu > \frac{2}{5}\chi ,$$
 or in other words $\chi > \frac{15}{4}\tau$. Reversing the 
 orientation of $M$, we  also have  $\chi > -\frac{15}{4}\tau$. 
 Since Lemma \ref{bshp} tells us that $\chi \leq 9$, we are therefore done.
 \end{proof}
 
 Our other main result now follows:
 
 \setcounter{main}{0} 
\begin{main}
Let $M$ be a smooth compact oriented 4-manifold with
(strictly) positive intersection form, and 
suppose that $g$ is an Einstein metric on $M$
which has non-negative sectional curvature. 
Then  $(M,g)$ is homothetically isometric
 to  $\bcp_{2}$, equipped with its standard Fubini-Study metric. 
\end{main}

\begin{proof}
By assumption,
 $b_+>0$ and $b_-=0$, so $\chi (M) = 2 + b_{+}$ and $\tau (M)= 
 b_{+} > 0$.
Hence
 $$\frac{15}{4} \tau = \frac{15}{4} b_{+} \geq \frac{9}{4} + b_{+} >
 2+ b_{+} = \chi.$$
  Theorem \ref{B} therefore insists that our Einstein metric $g$ 
  of non-negative sectional  curvature must satisfy either 
  $W^{+}\equiv 0$ or $W^{-}\equiv 0$. But since $\tau > 0$, 
  the signature formula (\ref{sig}) thus forces  $W^{-}\equiv 0$ and 
  $W^{+}\not\equiv 0$. In particular, $g$ is not flat,
  and, since it has non-negative  
  sectional curvature, its scalar curvature must somewhere be positive. 
 Thus  $(M,g)$ is a non-conformally-flat,
  self-dual Einstein  4-manifold of positive scalar
 curvature. A celebrated result of 
 Hitchin
 \cite[Theorem 13.30]{bes}, 
 originally discovered via twistor methods \cite{hit,FK}, 
  therefore tells us  that  
 $(M,g)$ must, up to isometry,
  be $\bcp_{2}$, equipped with a constant multiple of the
   Fubini-Study metric. 
  \end{proof}
  
\section{Einstein Constants} 

Given a smooth compact $n$-manifold $M$, for what values of 
$\lambda $ do Einstein's equations (\ref{eins}) have 
a unit-volume solution? The collection of all such $\lambda$ is called
\cite{bes}  the
set of {\em Einstein constants} for $M$, and 
constitutes an interesting smooth invariant of the manifold.
Assuming that $n > 2$, 
this set is just the  collection of 
 critical values  of the Riemannian functional 
 ${\cal S}/n$, where 
$${\cal S}(g) = \frac{\int_{M}s_{g}d\mu_{g}}{\left(
\int_{M}d\mu_{g}\right)^{\frac{n-2}{n}}},$$
since a metric is a critical point of $\cal S$ iff it is Einstein.

For $\bcp_{2}$, it is a relatively recent result \cite{leb4,GL}
that the set of Einstein constants has a maximal element, 
represented by, and only by, the Fubini-Study metric. 
This provides one new explanation for the rigidity \cite{koiso}
of the Fubini-Study metric. Theorem \ref{A} of course provides 
an ostensibly different explanation of this phenomenon, 
since the positivity of sectional curvatures 
is an open condition in the $C^{2}$ topology.  
However, the proof of Theorem \ref{A} tells us more. 
While it does not rule out the existence 
of an Einstein metric on $\bcp_2$ with sectional curvatures
of varying sign, it does  assert
 that the maximal element in the set of Einstein constants for $\bcp_{2}$
 is 
 {\em isolated}.
Moreover,  the form of this assertion is actually quantitative.

 \setcounter{main}{2}
\begin{thm}\label{C}
Suppose that $\bcp_2$ admits an Einstein metric $g$
which is not isometric to any multiple of  the Fubini-Study metric
$g_{0}$. Then 
$${\cal S}(g) <   \frac{1}{\sqrt{3}}{\cal S}(g_{0}),$$
where ${\cal S}(g_{0})=  
12\pi\sqrt{2}$.
\end{thm}

\begin{proof}
We may assume that the scalar curvature of $g$ is positive,
since the result is otherwise trivial.  By Hitchin's theorem 
\cite{hit,leb4}, we therefore know that $W^{-}_{g}\not\equiv 0$. 
For $\bcp_2, \chi=3$ and $\tau=1$. Thus (\ref{ii})
tells us that
$$\frac{1}{4\pi^{2}}
 	\int_{\bcp_{2}} \frac{s^2_g}{24} d\mu_{g} \leq 
 	\frac{2\chi - 3\tau }{3} (\bcp_{2}) =1.  $$
 	Moreover, the inequality is strict, since 
 	Corollary \ref{gapkid} would otherwise predict that 
 	the universal cover of 
 	$(\bcp_{2}, g)$ is {\em reverse-oriented} K\"ahler,
 	contradicting the fact that 
 	$\bcp_{2}$ is simply connected and has $b_{-}=0$.
 	Thus 
 	$${\cal S}(g) = \left(\int_{\bcp_{2}} s^{2}_{g}d\mu_{g}\right)^{{1/2}}
 	< \sqrt{4\pi^{2}\cdot 24} = 4\pi \sqrt{6} .$$
Since ${\cal S}(g_{0})=  
12\pi\sqrt{2}$, the result follows.
\end{proof}

While Bishop's inequality immediately implies that the 
set of Einstein constants for $S^{4}$ has a maximal element, 
represented by, and only by, the ``round'' metric,
the $S^{4}$-analog of Theorem \ref{C} was   only recently
proved  \cite{G1}.  
Unfortunately, however, an $S^{4}$-analog of 
Theorem \ref{A} remains 
out of our reach.  The best our present techniques can offer 
in this direction is

\begin{thm}\label{D}
Let $g$ be 
an Einstein metric of
non-negative sectional curvature
on $S^4$. If $g$ is not isometric to some multiple of the standard
 metric $g_{1}$, then 
 $$\frac{1}{\sqrt{5}}{\cal S}(g_{1})
<{\cal S}(g) < 
\frac{1}{\sqrt{3}}{\cal S}(g_{1}),
$$
where ${\cal S}(g_{1})= 8\pi \sqrt{6}$. 
\end{thm}
 
\begin{proof}
If $g$ is Einstein but has non-constant curvature, it cannot be 
conformally flat. By  (\ref{ii}), we thus have  
\[
 \int_{S^4} s^{2}_{g} d\mu _{g} \leq 
  4\pi^{2} \cdot 24 \frac{2\chi - 3\tau}{3}(S^{4})= 2^{7}\pi^{2}
\]
and the 
inequality is in fact strict, since $S^{4}$ cannot admit 
a (reverse-oriented) K\"ahler metric. Thus 
$${\cal S}(g) = \left(\int_{S^4} s^{2}_{g} d\mu _{g}\right)^{1/2}
< 8\pi \sqrt{2}.$$ 
On the other hand, if   $g$ has non-negative sectional curvature, 
(\ref{euler}) tells us    that 
$$
\frac{5}{8\pi^{2}}\int \frac{s^2_g}{24} d\mu _g > \chi (S^{4}) =2 ,
$$
so that 
$${\cal S}(g) = \left(\int_{S^4} s^{2}_{g} d\mu _{g}\right)^{1/2}
>\sqrt{ \frac{2^{7}\cdot 3\pi^{2}}{5}} =8\pi\sqrt{\frac{6}{5}}. $$
Since ${\cal S}(g_{1})= 8\pi \sqrt{6}$, the assertion follows. 
\end{proof}

\end{document}